\documentclass[12pt]{amsart}
\usepackage{amscd}
\def\NZQ{\Bbb}               

\def\ZZ{{\NZQ Z}}

\def\FF{{\NZQ F}}
\def\GG{{\NZQ G}}
\def\HH{{\NZQ H}}
\def\frk{\frak}               

\def\mm{{\frk m}}

\def\Phi{{\frk n}}
\def\Phi{{\frk N}}
\def\opn#1#2{\def#1{\operatorname{#2}}} 
\opn\chara{char}
\opn\length{\ell}
\opn\pd{pd}
\opn\rk{rk}
\opn\projdim{proj\,dim}
\opn\injdim{inj\,dim}
\opn\rank{rank}
\opn\depth{depth}
\opn\grade{grade}
\opn\height{height}
\opn\embdim{emb\,dim}
\opn\codim{codim}

\opn\Tr{Tr}
\opn\bigrank{big\,rank}
\opn\superheight{superheight}\opn\lcm{lcm}
\opn\trdeg{tr\,deg}%
\opn\reg{reg}
\opn\lreg{lreg}
\opn\ini{in}
\opn\lpd{lpd}
\opn\size{size}
\opn\div{div}
\opn\Div{Div}
\opn\cl{cl}
\opn\Cl{Cl}
\opn\Spec{Spec}
\opn\Supp{Supp}
\opn\supp{supp}
\opn\Sing{Sing}
\opn\Ass{Ass}
\opn\Ann{Ann}
\opn\Rad{Rad}
\opn\Soc{Soc}
\opn\Im{Im}
\opn\Ker{Ker}
\opn\Coker{Coker}
\opn\Am{Am}
\opn\Hom{Hom}
\opn\Tor{Tor}
\opn\Ext{Ext}
\opn\End{End}
\opn\Aut{Aut}
\opn\id{id}

\opn\nat{nat}
\opn\pff{pf}
\opn\Pf{Pf}
\opn\GL{GL}
\opn\SL{SL}
\opn\mod{mod}
\opn\ord{ord}
\opn\Gin{Gin}
\opn\Hilb{Hilb}
\opn\aff{aff}
\opn\con{conv}
\opn\relint{relint}
\opn\st{st}
\opn\lk{lk}
\opn\cn{cn}
\opn\core{core}
\opn\vol{vol}
\opn\link{link}
\opn\star{star}
\opn\gr{gr}

\def\pot#1#2{#1[\kern-0.28ex[#2]\kern-0.28ex]}

\opn\dirlim{\underrightarrow{\lim}}
\opn\inivlim{\underleftarrow{\lim}}
\let\union=\cup
\let\sect=\cap
\let\dirsum=\oplus

\let\iso=\cong
\let\Union=\bigcup
\let\Sect=\bigcap

\let\to=\rightarrow
\let\To=\longrightarrow
\def\Implies{\ifmmode\Longrightarrow \else
        \unskip${}\Longrightarrow{}$\ignorespaces\fi}
\def\implies{\ifmmode\Rightarrow \else
        \unskip${}\Rightarrow{}$\ignorespaces\fi}
\def\iff{\ifmmode\Longleftrightarrow \else
        \unskip${}\Longleftrightarrow{}$\ignorespaces\fi}

\let\:=\colon
\newtheorem{Theorem}{Theorem}[section]
\newtheorem{Lemma}[Theorem]{Lemma}
\newtheorem{Corollary}[Theorem]{Corollary}

\newtheorem{Remark}[Theorem]{Remark}

\newtheorem{Definition}[Theorem]{Definition}

\let\epsilon\varepsilon
\let\phi=\varphi
\let\kappa=\varkappa
\def\qed{\ifhmode\textqed\fi
      \ifmmode\ifinner\quad\qedsymbol\else\dispqed\fi\fi}
\def\textqed{\unskip\nobreak\penalty50
       \hskip2em\hbox{}\nobreak\hfil\qedsymbol
       \parfillskip=0pt \finalhyphendemerits=0}
\def\dispqed{\rlap{\qquad\qedsymbol}}

\opn\dis{dis}
\def\pnt{{\raise0.5mm\hbox{\large\bf.}}}

\opn\Lex{Lex}

\def\Coh#1#2{H_{\mm}^{#1}(#2)}
\def\hchst#1{for all $u\in G(I)$ there  exists $i\notin {#1}$
such that $\nu_i(u) > a_i \geq 0$}

\begin{document}

\title{On the radical of a monomial ideal}

\author{J\"urgen Herzog,  Yukihide Takayama and Naoki Terai}
\address{J\"urgen Herzog, Fachbereich Mathematik und
Informatik, Universit\"at Duisburg-Essen, Campus Essen,
45117 Essen, Germany}
\email{juergen.herzog@uni-essen.de}

\address{Yukihide Takayama, Department of Mathematical
Sciences, Ritsumeikan University,
1-1-1 Nojihigashi, Kusatsu, Shiga 525-8577, Japan}
\email{takayama@se.ritsumei.ac.jp}

\address{Naoki Terai,Department of Mathematics, Faculty of Culture and Education, Saga University, Saga 840-8502, Japan }
\email{terai@cc.saga-u.ac.jp}
\maketitle

\begin{abstract}
Algebraic and combinatorial properties of a monomial ideal and its radical are compared.
\end{abstract}

\section{Introduction}
There are simple examples of Cohen-Macaulay ideals whose radical is not Cohen-Macaulay. The first such example is probably due to
Hartshorne \cite{H}, who proved that in positive characteristic the toric ring
$K[s^4,s^3t,st^3,t^4]$ is a set theoretic complete intersection. With CoCoA or other computer algebra systems many other examples,
 also in characteristic zero,  can be constructed. The following example due Conca was computed with  CoCoA:  let $S=K[x_1,x_2,x_3,x_4, x_5]$
 and  $J=(x_2^2 - x_4x_5, x_1x_3 -  x_3x_4, x_3x_4 - x_1x_5)\subset S$. Then   $S/J$ is a 2-dimensional Cohen-Macaulay ring,
 $\sqrt{J}= (x_1x_3 - x_1x_5,  x_3x_4 - x_1x_5,  x_2^2 -
x_4x_5,  x_1^2x_2 - x_1x_2x_4,  x_2x_3^2 - x_2x_3x_5)$ and $S/\sqrt{J}$ is {\em not} Cohen-Macaulay. Indeed, the  depth of $S/\sqrt{J}$ equals $1$. On the other hand it is well-known that the Cohen-Macaulay property of a monomial ideal is inherited by its radical.  The reason is that the radical of a monomial ideal is essentially obtained by polarization and localization. This observation,  was communicated to the third author by David Eisenbud. Both operations, polarization and localization, preserve the Cohen-Macaulay property. An explicit proof of this fact can be found in \cite{Ta}. The purpose of this paper is to exploit this idea and to show that many other nice properties are inherited by the radical of a monomial ideal.

\section{The comparison}

For the proof of the main result of this paper we need  some preparation. We begin with the following extension \cite[Theorem 1.1]{T} of
Hochster's formula \cite[Theorem 5.3.8]{BH} describing  the local cohomology of a monomial ideal.

Let $K$ be a field, $S=K[x_1,\ldots, x_n]$ the polynomial ring and $I\subset S$ a monomial ideal. The unique minimal monomial system of
generators of $I$ is denoted by $G(I)$. For $i=1,\ldots, n$ we set
\[
t_i=\max\{\nu_i(u)\: u\in G(I)\},
\]
where for a monomial $u\in S$, $u=x_1^{a_1}\cdots x_n^{a_n}$ we set $\nu_i(u)=a_i$ for $i=1,\ldots,n$.

For $a=(a_1,\ldots, a_n)\in \ZZ^n$, we set
\[
G_a= \{i \: 1\leq i\leq n,\; a_i < 0\},
\]
and define the simplicial complex $\Delta_a(I)$ whose faces are the sets $L\setminus G_a$ with $G_a\subset L$, and such that $L$ satisfies the
following condition: {\hchst{L}.

Notice that the inequality $a_i\geq 0$
in the definition of $\Delta_a(I)$
follows from the condition $i\notin L\supset G_a$. It  is  included only for the reader's convenience.

\medskip
\noindent
With the notation introduced one has

\begin{Theorem}[Takayama \cite{T}]
\label{hochster}
Let $I\subset S$ be a monomial ideal. Then
the Hilbert series of the local cohomology modules of $S/I$
with respect to the $\ZZ^n$-grading is given by
\begin{equation*}
{\Hilb}(\Coh{i}{S/I}, {\bf t})
=\sum_{F\in\Delta}
 \sum_{a}
   \dim_K\tilde{H}_{i-\vert F\vert -1}(\Delta_a(I); K) {\bf t}^a
\end{equation*}
where $\Delta$ is the
simplicial complex corresponding to the Stanley-Reisner
ideal $\sqrt{I}$, and the second sum is taken over all $a\in\ZZ^n$ such that $a_i\leq t_i-1$ for all $i$, and $G_a=F$.
\end{Theorem}

As a first application of this theorem we have

\begin{Corollary}
\label{bound}
Let $I\subset S$ be a monomial ideal. Then
\begin{equation*}
a(S/I) \leq \sum_{i=1}^n t_i - n,
\end{equation*}
where $a(S/I)$ is the $a$-invariant  of $S/I$.
\end{Corollary}

\begin{proof}
By Theorem~\ref{hochster}, we know that $\Coh{i}{R}_{a} =0$ for
all $i$ and for all $a\in\ZZ^n$ such that $a_i>t_i-1$ for some
$i$. Thus in particular, if $d=\dim R$, then $\Coh{d}{R}_{j}=0$
for $j
>  \sum_{i=1}^n t_i -n$.
\end{proof}

We say that $S/I$ has {\em maximal $a$-invariant} if the upper bound in Corollary \ref{bound} is attained, that is, if
$a(S/I)=\sum_{i=1}^n t_i -n$.

\medskip
For our main theorem the next corollary is important.

\begin{Corollary}
\label{useful}
Let $I \subset S$ be a monomial ideal. Then
we have the following isomorphisms of $K$-vector spaces
\begin{equation*}
\Coh{i}{S/I}_{a} \iso \Coh{i}{S/\sqrt{I}}_{a}
\end{equation*}
for all $a\in \ZZ^n$ with $a_i \leq 0$ for $1\leq i\leq n$.
\end{Corollary}

\begin{proof}
Consider the multigraded Hilbert series of $\Coh{i}{S/I}$
and $\Coh{i}{S/\sqrt{I}}$. Let $a\in\ZZ^n$ be such that
$a_i \leq 0$ for all $1\leq i\leq n$.  Then
by Theorem~\ref{hochster}, we have
\begin{eqnarray*}
\dim_K\Coh{i}{S/I}_a
   & = & \dim_K\tilde{H}_{i-\vert F\vert -1}(\Delta_a(I); K), \quad \text{and}   \\
\dim_K\Coh{i}{S/\sqrt{I}}_a
   & = & \dim_K\tilde{H}_{i-\vert F\vert -1}(\Delta_a(\sqrt{I}); K),
\end{eqnarray*}
For a monomial $u$ we set $\supp(u)=\{i\: x_i \text{ divides } u\}$. Now since for every $u\in G(I)$ there exists $v\in G(\sqrt{I})$ such
that $\supp(u) \supset \supp(v)$, and
since
for every $v\in G(\sqrt{I})$ there exists $u\in G(I)$ such that
$\supp(v) = \supp(u)$, it follows  that
$\Delta_a(I) = \Delta_a(\sqrt{I})$. Thus we have
$\dim_K \Coh{i}{S/I}_a = \dim_K \Coh{i}{S/\sqrt{I}}_a$.
\end{proof}

Let $M$ be a graded $S$-module. For the convenience of the reader we recall the following two concepts which generalize the Cohen-Macaulay property and non-pure shellability  of simplicial complexes.

The following definition is due to Stanley \cite[Section II, 3.9]{St}:

\begin{Definition}
\label{stanley}
{\em Let $M$ be a finitely generated graded $S$-module. The module $M$ is {\em sequentially
Cohen-Macaulay} if there exists a finite filtration
\[
0=M_0\subset M_1\subset M_2\subset\ldots\subset M_r=M
\]
of $M$ by graded submodules of $M$ such that each quotient $M_i/M_{i-1}$ is CM, and  $\dim
M_1/M_0<\dim M_2/M_1<\ldots <\dim M_r/M_{r-1}$.}
\end{Definition}

It is known (see for example \cite[Corollary 1.7]{HS}) that if $M$ is sequentially Cohen-Macaulay, then the filtration given in the definition is uniquely determined. We call it the {\em attached filtration} of the sequentially Cohen-Macaulay module $M$.

The uniqueness of the filtration is seen as follows: suppose $\depth M=t$, then $M_1$ is the image of the natural map $\Ext_S^{n-t}(\Ext_S^{n-t}(M,\omega_S),\omega_S)\to M$. Here $\omega_S=S(-n)$ is the canonical module of $S$. Then one notices that $M/M_1$ is again sequentially Cohen-Macaulay and uses induction on the length of the attached sequence.

In case $M$ is a cyclic module, say, $M=S/I$, with attached filtration $0=M_0\subset M_1\subset \cdots$, each of the the modules $M_i$ is an ideal in $S/I$, and hence is of the form $I_i/I$ for certain (uniquely determined) ideals $I_i\subset S$. Thus $S/I$ is sequentially Cohen-Macaulay, if and only of there exists a chain of graded ideals
\[
I=I_0\subset I_1\subset I_2\subset \ldots \subset I_r=S
\]
such that each factor module $I_{i+1}/I_i$ is Cohen-Macaulay with $$\dim I_{i+1}/I_i<\dim I_{i+2}/I_{i+1}$$ for $i=0,\ldots, r-2$. Moreover if this property is satisfied, then this chain of ideals is uniquely determined.

In the particular case that $I$ is a monomial ideal, the natural map
$$\Ext_S^{n-t}(\Ext_S^{n-t}(S/I,\omega_S),\omega_S)\to S/I$$
is a homomorphism of multigraded $S$-modules. This implies that the attached chain of ideals of the
sequentially Cohen-Macaulay module $S/I$  is a chain of monomial ideals.

\medskip
Now let us briefly describe the other concept which was introduced by Dress \cite{D}:

\begin{Definition}
\label{Dress}
{\em Let $M$ be a finitely generated graded $S$-module. A filtration
\[
0=M_0\subset M_1\subset M_2\subset\ldots\subset M_r=M
\]
of $M$ by graded submodules of $M$ is called {\em clean} if for all $i=1,\ldots, r$ there exists a minimal prime ideal $P_i$ of $M$  such that $M_i/M_{i-1}\iso S/P_i$. The module $M$ is called {\em clean} if  it has a clean filtration.}
\end{Definition}

Again, if $M=S/I$ is cyclic, then $S/I$ is clean if there exists a chain of ideals $I=I_0\subset I_1\subset I_2\subset\ldots\subset I_{r-1}\subset I_r=S$ such that $I_{i+1}/I_i\iso S/P_i$ with $P_i$ a minimal prime ideal of $I$. In other words, for all $i=0,\ldots, r-1$ there exists $f_{i+1}\in I_{i+1}$ such that $I_{i+1}=(I_i,f_{i+1})$ and $P_i=I_i:f_{i+1}$. In case $I$ is a monomial ideal we require that all $f_i$ are monomials.

Dress \cite{D} shows that a Stanley-Reisner ideal $I_\Delta$ is
clean if and only the simplicial complex $\Delta$ is non-pure
shellable in the sense of Bj\"orner and Wachs \cite{BW}.

\medskip
In the proof of our main theorem we use polarization, as indicated
in the introduction. Let $I=(u_1,\ldots, u_m)$ with
$u_i=x_1^{a_{i1}}\cdots x_n^{a_{in}}$. We fix some number $i$ with
$1\leq i\leq n$,  introduce a new variable $y$, and set
$v_k=x_1^{a_{k1}}\cdots x_i^{a_{ki}-1}y\cdots x_n^{a_{kn}}$ if
$a_{ki}>1$,and $v_k=u_k$ otherwise. We call  $J=(v_1,\ldots, v_m)$
the {\em $1$-step polarization of $I$ with respect to the variable
$x_i$}. The element $y-x_i$ is regular on $S[y]/J$ and
$(S[y]/J)/(y-x_i)(S[y]/J)\iso S/I$, see \cite[Lemma 4.2.16]{BH}.

Let as above $t_i=\max\{\nu_i(u_j)\: j=1,\ldots,m\}$, and set
$t=\sum_{i=1}^nt_i-n$. Then it is clear that if we apply $t$
suitable 1-step polarizations, we end up with a squarefree
monomial ideal $I^p$, which is called the {\em complete
polarization of $I$}.

\medskip
Now we are ready to present the  main result of this section.

\begin{Theorem}
\label{main}
Let $K$ be a field, $S=K[x_1,\ldots,x_n]$ the polynomial ring over $K$, and $I\subset S$ a monomial ideal. Suppose that   $S/I$ satisfies one of the following properties: $S/I$ is {\em (i)} Cohen-Macaulay,
{\em (ii)} Gorenstein, {\em (iii)} sequentially Cohen-Macaulay, {\em (iv)} generalized Cohen-Macaulay, {\em (v)} Buchsbaum, {\em (vi)} clean,
or {\em (vii)} level and has maximal $a$-invariant.
Then $S/\sqrt{I}$ satisfies the corresponding property.
\end{Theorem}

\begin{proof}
We  first use the trick, mentioned in the introduction, to show that the  Betti-numbers $\beta_{i}(I)$ of $I$ do not increase when passing to $\sqrt{I}$.

We denote by $I^p$ the complete polarization of $I$. Let $T$ be the polynomial ring in the variables that are needed to polarize $I$. Then $I^p$ is a squarefree monomial ideal in $T$ with $\beta_{i}(I^p)=\beta_{i}(I)$ for all $i$. It is easy to see that if we localize at the multiplicative set $N$ generated by the  {\em new} variables which are needed to polarize $I$, one obtains $I^pT_N=(\sqrt{I})T_N$. Since localization is an exact functor, the localized free resolution will be a possibly non-minimal free resolution of $(\sqrt{I})T_N$. Since the extension $S\to T_N$ is flat, the desired inequality follows.

\medskip
Proof of (i) and (ii): The inequality $\beta_{i}(\sqrt{I})\leq \beta_i(I)$ implies that $\depth S/\sqrt{I}\geq \depth S/I$. On the other hand, $\dim S/I=\dim S/\sqrt{I}$. This implies that $S/\sqrt{I}$ is Cohen-Macaulay, if $S/I$ is so.

Suppose now that $S/I$ is Gorenstein. Then $\beta_q(S/I)=1$ where
$q$ is the codimension of $I$, see \cite[Theorem 3.3.7 and
Corollary 3.3.9]{BH}. Therefore, $\beta_q(S/\sqrt{I})\leq 1$.
Since $I$ and $\sqrt{I}$ have the same codimension, we see that
$\beta_q(S/\sqrt{I})>0$, and hence $\beta_q(S/\sqrt{I})=1$. Again
using \cite[Theorem 3.3.7 and Corollary 3.3.9]{BH} we conclude
that $S/\sqrt{I}$ is Gorenstein. This fact follows also from
\cite[Corollary 3.4]{BH1}.

\medskip
Proof of (iii): Since $S/I$ is sequentially Cohen-Macaulay there exists a chain of monomial ideals
\[
I=I_0\subset I_1\subset I_2\subset \cdots \subset I_k=S
\]
such that $I_{j+1}/I_j$ is Cohen-Macaulay for all $j=0,\ldots,k-1$ and such that $\dim I_1/I_0< \dim I_2/I_1< \ldots <\dim I_k/I_{k-1}$.

Suppose $x_1^a$ with $a>1$ divides a generator of $I$. Then we apply a $1$-step polarization for $x_1$ to all the ideals $I_i$, and obtain a chain of ideals $J=J_0\subset J_1\subset J_2\subset \cdots \subset J_k =\tilde{S}$ where $\tilde{S}=S[y]$. It follows that $y-x_1$ is $\tilde{S}/J_i$-regular and $(\tilde{S}/J_i)/(y-x_1)(\tilde{S}/J_i)\iso S/I_i$ for all $i$. Therefore $y-x_1$ is $J_{i+1}/J_i$-regular, and  $(J_{i+1}/J_i)/(y-x_1)(J_{i+1}/J_i)\iso I_{i+1}/I_i$. Thus $J$ is  sequentially Cohen-Macaulay.

Since the complete polarization $I^p_i$ of the ideals $I_i$ for $i=1,\ldots k$, is obtained by a sequence of 1-step polarizations, it follows that $I^p$ is sequentially Cohen-Macaulay. As $I_i^p/I_{i+1}^p$ is Cohen-Macaulay, we conclude as in the proof of (i) that $\sqrt{I_{i+1}}/\sqrt{I_i}$ is Cohen-Macaulay of the same dimension as  $I_{i+1}/I_i$. This shows that $\sqrt{I}$ is sequentially Cohen-Macaulay.

\medskip
Proof of (iv) and (v): Assuming that $S/I$ is generalized Cohen-Macaulay or Buchsbaum, one has  that $S/I$ is equidimensional and that $H^i_\mm(S/I)_j=0$ for all $i<\dim S/I$, and all but finitely many $j$. Since $I$ and $\sqrt{I}$ have the same minimal prime ideals, it follows that $\sqrt{I}$ is again equidimensional.

Let $\ZZ^n_-$ be the set of all  $a\in \ZZ^n$ such that $a_i\leq 0$ for $i=1,\ldots,n$.
By Corollary \ref{useful}, $H^i_\mm(S/I)_a=H^i_\mm(S/\sqrt{I})_a$ for all $a\in \ZZ^n_-$. Moreover, by Hochster's formula, $H^i_\mm(S/\sqrt{I})_a=0$ for all $a\not\in  \ZZ^n_-$.
Therefore, $\dim_K H^i_\mm(S/\sqrt{I})_j\leq \dim_K H^i_\mm(S/I)_j$ for all $j\leq 0$ and $H^i_\mm(S/J)_j=0$ for $j>0$. It is known \cite{Sch}  that a squarefree monomial ideal is Buchsbaum if and only if it is generalized Cohen-Macaulay. Thus (iv) and (v) follow.

\medskip
Proof of (vi): Assuming that $S/I$ is clean, there exists a chain of monomial ideals $I=I_0\subset I_1\subset I_2\subset\ldots\subset I_{r-1}\subset I_r=S$ such that $I_{i+1}/I_i\iso S/P_i$ with $P_i$ a minimal prime ideal of $I$. We claim that $\sqrt{I_{i+1}}/\sqrt{I_i}=S/P_i$, if $\sqrt{I_{i+1}}\neq \sqrt{I_{i}}$. This then implies that $S/\sqrt{I}$ is clean, since the prime ideals  $P_i$ are also minimal prime ideals of $\sqrt{I}$.

In order to prove this claim we introduce some notation: let $u=x_1^{a_1}x_2^{a_2}\cdots x_n^{a_n}$ and $v= x_1^{b_1}x_2^{b_2}\cdots x_n^{b_n}$ be two monomials. Then we set
\[
u:v=\prod_{i=1}^nx_i^{\max\{a_i-b_i,0\}},\quad\text{and}\quad u_{red}=\prod_{i\atop a_i>0}x_i.
\]
We then have
\begin{eqnarray}
(u:v)_{red}= (u_{red}:v_{red})\prod_{i,\atop a_i>b_i>0}x_i.
\end{eqnarray}
Note that if $I$ is a monomial ideal with monomial generators $u_1,\ldots, u_m$, then
\[
\sqrt{I}=((u_1)_{red},\ldots, (u_m)_{red})\quad \text{and}\quad I:v=(u_1:v, \ldots, u_m:v).
\]
Back to the proof of our claim, our assumption implies that for  all $i=0,\ldots, r-1$ there exists a monomial $v_{i+1}\in I_{i+1}$ such that $I_{i+1}=(I_i,v_{i+1})$ and $P_i=I_i:v_{i+1}$. Suppose $P_i=(x_{i_1},\ldots, x_{i_s})$. Then  $P_i=I_i:v_{i+1}$ if and only if
\begin{enumerate}
\item[(a)] for all $j=1,\ldots,s$ there exists $u\in I_i$ such that $u:v_{i+1}=x_{i_j}$, and
\item[(b)] for all monomial generators $w\in I_i$ there exists an integer $j$ with $1\leq j\leq s$ such that $x_{i_j}|(w:v_{i+1})$.
\end{enumerate}

We need to show that $P_i=\sqrt{I_i}:(v_{i+1})_{red}$, if $(v_{i+1})_{red}\not\in \sqrt{I_i}$, and prove this by checking (a) and (b) for the pair $\sqrt{I_i}$ and $(v_{i+1})_{red}$.

Let $j$ be an integer with $1\leq j\leq s$. Then there exists $u\in I_i$ such that $u:v_{i+1}=x_{i_j}$.
Suppose $u=\prod_{k=1}^nx_k^{a_k}$ and $v_{i+1}= \prod_{k=1}^nx_k^{b_k}$, then
 (1) implies that $x_{i_j}=(u:v_{i+1})_{red}=(u_{red}:(v_{i+1})_{red})w$ where $w=\prod_{k,\; a_k>b_k>0}x_k$.
Suppose  $x_{i_j}$ divides $w$, then  $u_{red}:(v_{i+1})_{red}=1$. This implies that $(v_{i+1})_{red}\in \sqrt{I_i}$, a contradiction. Therefore $u_{red}:(v_{i+1})_{red}=x_{i_j}$, and this proves (a). The argument also shows that $b_{i_j}=0$ for $j=1,\ldots, s$.

For the proof of (b), let $w\in I_i$ be a monomial generator. Then there exists an integer   $j$ with $1\leq j\leq s$ such that $x_{i_j}|(w:v_{i+1})$. It follows that $x_{i_j}$ divides $(w:v_{i+1})_{red}$. Let $w=\prod_{k=1}^n x_k^{c_k}$. Then (1) implies that $x_{i_j}$ divides $(w_{red}:(v_{i+1})_{red})\prod_{k,\; c_k>b_k>0}x_k$. However, $b_{i_j}=0$, as we have seen in the proof of (a). Therefore, $x_{i_j}$ divides $(w_{red}:(v_{i+1})_{red})$. Since $\sqrt{I_i}$ is generated by the monomials $w_{red}$ where the monomials $w$ are the generators of $I_i$, condition (b) follows.

\medskip
Proof of (vii): By assumption $S/I$ is level. This means  that
$S/I$ is Cohen-Macaulay and that all generators of the canonical
module $\omega_{S/I}$ of $S/I$ have the same degree, say $g$. In
this situation the $a$-invariant $a(S/I)$ of $S/I$ is just $-g$,
see \cite[Section 3.6]{BH}. Suppose $d= \dim S/I$; then $I$ has a
graded minimal free resolution $\FF$ of length $q=n-d-1$ with
$F_q=S^b(-c)$, Since $\omega_{S/I}$ may be represented as the
cokernel of $F_{q-1}^*\to F_q^*$, which is dual  of the map
$F_q\to F_{q-1}$ with respect to $S(-n)$, it follows that
$a(S/I)=c-n$.

For $i=1,\ldots, n$ we set again
\[
t_i=\max\{\nu_i(u)\: u\in G(I)\}.
\]
By Corollary \ref{bound},  one has  the upper bound $a(S/I)\leq
\sum_{i=1}^n t_i - n$. Since we assume that $S/I$ has maximal
$a$-invariant, the upper bound is reached.  Let $I^p\subset T$ the
complete polarization of $I$. This polarization requires precisely
$t= \sum_{i=1}^n t_i - n$ 1-step polarizations. It follows that
$S/I$ is obtained from $T/I^p$ as a residue class ring modulo a
regular sequence of linear forms of length $t$. From the above
description of the $a$-invariant we now conclude   that $a(T/I
^p)=a(S/I)-t=0$. Let $\GG$ be the multigraded minimal free
resolution of the squarefree monomial ideal $I^p$. Since $\projdim
I^p=\projdim I=q$, and since $a(T/I^p)=0$, we see that
$G_q=T(-m)^b$, where $m=n+t=\dim T$. This implies that $G_q$ as a
multigraded module is isomorphic to $T(-e)^b$ where
$e=(1,1,\ldots,1)$.

For $i=1,\ldots,m$ let $e_i$ be the $i$th canonical basis vector
of $\ZZ^m$. Then $e=\sum_{i=1}^me_i$, and we may assume that $\deg
x_i=e_i$ for $i=1,\ldots,n$, while the new variables have the
multidegrees $e_i$ with $i=n+1,\ldots, m$.  We define a new
multigrading on $T$ and $T/I^p$: for an element $f$ of multidegree
$a$ we set $\deg' f=\pi(a)$, where $\pi\: \ZZ^m\to \ZZ^n$ is the
projection onto the first $n$ components of $\ZZ^m$.

As above, let $N$ be the multiplicative set generated by the $t$ new variables which are needed to polarize $I$.
Then $I^pT_N=\sqrt{I}T_N$, and localization with respect to $N$ preserves the new multigrading since $\deg' f=0$
for all $f\in N$. Therefore $\GG_N$ is, with respect to the new grading,  a multigraded free $T_N$-resolution of
$\sqrt{I}T_N$ with $(G_q)_N=T_N(-1,\ldots, -1)^b$ and $(-1,\ldots, -1)\in\ZZ^n$.

Let $\HH$ be the multigraded minimal free $S$-resolution of
$\sqrt{I}$.  Then $\HH T_N$ is the minimal multigraded free
$T_N$-resolution of $\sqrt{I}T_N$.  A comparison with the
(possibly non-minimal) graded free $T_N$-resolution $\GG_N$ shows
that $H_q$ is a direct summand of copies of $S(-1,\ldots,-1)$.
Since $S/I$ and $S/\sqrt{I}$ are Cohen-Macaulay of the same
dimension, we see that $q=\projdim I=\projdim \sqrt{I}$. Therefore
all summands in the last step of the resolution $\HH$ of
$S/\sqrt{I}$ have the same shift. This show that $S/\sqrt{I}$ is
level.
\end{proof}

\begin{Remark}
{\em In Theorem 2.6(i) (or (iv)), it suffices to require that $I$ is an arbitrary homogeneous (generalized) Cohen-Macaulay ideal whose radical $\sqrt I$ is a monomial
ideal, i.e.\ we do not need to require that $I$ itself is a monomial ideal.

Indeed  it is enough  to prove that there is a surjective homomorphism $H^i_\mm(S/I) \longrightarrow  H^i_\mm(S/\sqrt{I})$ for all $i$. The natural surjective
map $S/I \longrightarrow S/ \sqrt I$ induce for all $i$   commutative diagrams
\[
\begin{array}{ccc}
\Ext^{i}(S/\sqrt I, S)  & \longrightarrow  & \Ext^{i}(S/I, S) \\
\downarrow  &  & \downarrow  \\
H_{\sqrt I}^{i}(S)        &  \longrightarrow &    H_{I}^{i}(S).
\end{array}
\]
Since $H_{\sqrt I}^{i}(S) \cong   H_{I}^{i}(S)$ and since $\Ext^{i}(S/\sqrt I, S) \longrightarrow  H_{\sqrt I}^{i}(S)$ is an essential extension (see
\cite{Te}), it follows that $\Ext^{i}_S(S/\sqrt I, S) \longrightarrow \Ext^i_S(S/I, S)$ is injective for all $i$. Hence the desired conclusion follows by
local duality.

On the other hand, as for the Gorenstein property, we must assume that $I$ is a monomial ideal. For example, $I=(xy+yz, xz)$ is a complete intersection, hence,
a Gorenstein ideal, while $\sqrt I = (xy, yz, xz)$  is not Gorenstein.}
\end{Remark}

\section{The inverse problem}

The results of the previous section indicate the following question: for a subset $F\subset [n]$, let $P_F$ be the prime ideal generated by the $x_i$ with $i\in F$. The minimal prime ideals of a squarefree $I$ are all of this form, and since $I$ is a radical ideal it is the intersection of its minimal prime ideals, say,   $I=\Sect_{i=1}^r P_{F_i}$ with $F_i\subset [n]$.

Suppose $I$ is Cohen-Macaulay. For which exponents $a_{ij}$ is the ideal
\[
J=\Sect_{i=1}^r (x_j^{a_{ij}}\: j\in F_i)
\]
again Cohen-Macaulay?

Of course if we raise the $x_i$ uniformly to some power, say $x_i$ is replaced by  $x_i^{a_i}$ everywhere in the intersection, then the resulting ideal $J$ is the image of the flat map $S\to S$ with $x_i\mapsto x_i^{a_i}$ for all $i$. Thus in this case $J$ will be Cohen-Macaulay, if $I$ is so. On the other hand, if we allow arbitrary exponents, the question seems to be quite delicate, and we do not know a general answer. However, if we require that for {\em all} choices of exponents  the resulting ideal is again Cohen-Macaulay, a complete answer is possible.

We need a definition to state the next result. Let $L$ be a monomial ideal. Lyubeznik \cite{L} defines the {\em size} of $L$ as follows: let $L=\Sect_{j=1}^rQ_j$ be an irredundant primary decomposition of $L$, where the $Q_i$ are monomial ideals. Let $h$ be the height  of $\sum_{j=1}^rQ_j$, and denote by $v$ the minimum number $t$ such that there exist $j_1,\ldots,j_t$ with $\sqrt{\sum_{i=1}^tQ_{j_i}}=\sqrt{\sum_{j=1}^rQ_j}$. Then $\size L=v+(n-h)-1$.

Since for monomial ideals the operations of forming sums and taking radicals  can be exchanged, the numbers $v$ and $h$, and hence the size of $L$ depends only on the associated prime ideals of $L$.

We shall need the following result of Lyubeznik \cite[Proposition 2]{L}:

\begin{Lemma}
\label{Lyubeznik}
Let $L$ be a monomial ideal in $S$. Then $\depth S/L\geq \size L$.
\end{Lemma}

Now we can state the main result of this section.

\begin{Theorem}
\label{converse}
Let $I\subset S=K[x_1,\ldots, x_n]$ be a  Cohen-Macaulay squarefree monomial ideal, and write
\[
I=\Sect_{i=1}^r P_{F_i},
\]
where the sets  $F_i\subset [n]$ are pairwise distinct, and  all have the same cardinality $c$.

For $i=1,\ldots, r$ and $j=1,\ldots, c$ we choose integers $a_{ij}\geq 1$, and set
 $$Q_{F_i}=(x_j^{a_{ij}}\: j\in F_i)\quad \text{for}\quad  i=1,\ldots, r.$$

Then the following conditions are equivalent:
\begin{enumerate}
\item[(a)] for all choices of the integers $a_{ij}$ the ideal
\[
J= \Sect_{i=1}^rQ_{F_i}
\]
is Cohen-Macaulay; \item[(b)] for each subset $A\subset [r]$, the
ideal $I_A=\Sect_{i\in A}P_{F_i}$ is Cohen-Macaulay; \item[(c)]
$\height P_{F_i}+P_{F_j}=c+1$ for all $i\neq j$; \item[(d)] for
$r\geq 2$ either $|\Union_{i=1}^r F_i|=c+1$, or $|\Sect_{i=1}^r
F_i|=c-1$; \item[(e)] after a suitable permutation of the elements
of $[n]$ we either have
$$F_i=\{1,\ldots,i-1,i+1,\ldots,c,c+1\}\quad\text{for}\quad
i=1,\ldots,r,$$ or
$$F_i=\{1,\ldots,c-1,c-1+i\}\quad \text{for}\quad  i=1,\ldots, r;$$
\item[(f)] $\size I=\dim S/I$;
\item[(g)] $S/L$ is Cohen-Macaulay for any monomial ideal $L$ such that $\Ass L= \Ass I$.
\end{enumerate}
 \end{Theorem}

\begin{proof}
(a)\implies (b): Let $Q_{F_i}=(x_j^2\: j\in F_i)$ if $i\in A$, and
$Q_{F_i}=P_{F_i}$ if $i\not\in A$. By assumption,
$J=\Sect_{i=1}^rQ_{F_i}$ is Cohen-Macaulay. Hence the complete
polarization $J^p$ of $J$ is again Cohen-Macaulay. We have $J^p=
\Sect_{i=1}^rQ_{F_i}^p$ with $Q_{F_i}^p=(x_jy_j\: j\in F_i)$ if
$i\in A$, and $Q_{F_i}^p=P_{F_i}$ if $i\not\in A$. Let $N$ be the
multiplicative set generated by all the variables $x_i$.  Then
$J^p_N$ is Cohen-Macaulay, and  hence
\[
J^p_N=\Sect_{i\in A}(y_j\: j\in F_i).
\]
This shows that $I_A=\Sect_{i\in A}P_{F_i}$ is Cohen-Macaulay.

(b)\implies (c): Consider the exact sequence
\[
0\To S/(P_{F_i}\sect P_{F_j})\To S/P_{F_i}\dirsum S/P_{F_j}\To S/(P_{F_i}+P_{F_j})\To 0.
\]
The rings $S/P_{_i}$ and $S/P_{_j}$ are Cohen-Macaulay of dimension $n-c$, while $S/(P_{F_i}+P_{F_j})$ is Cohen-Macaulay of dimension $n-d$ where $d$ is the height of $P_{F_i}+P_{F_j}$. The exact sequence yields that $S/(P_{F_i}\sect P_{F_j})$ is Cohen-Macaulay if and only if $d=c+1$.

Since by assumption $S/P_{F_i}\sect P_{F_j}$ is Cohen-Macaulay for all $i\neq j$, the assertion follows.

(c)\implies (d): We must show: given a collection of subsets $F_1,\ldots, F_r\subset [n]$ with
\begin{enumerate}
\item[(i)] $|F_i|=c$ for all $i$;
\item[(ii)] $|F_i\union F_j|=c+1$ for all $i\neq j$.
\end{enumerate}
Then either  $|\Union_{i=1}^rF_i|=c+1$, or $|\Sect_{i=1}^r F_i|=c-1$.

Suppose this is not the case. Then, since $|F_1\sect F_2|=c-1$ and $|F_1\union F_2|=c+1$, there exist  integers $i$ and $j$ such that $F_1\sect F_2\not\subset F_i$,  and $F_j\not\subset F_1\union F_2$.
The conditions (i) and (ii) then imply that there exists an element $x\in F_1\sect F_2$ such that $F_1\union F_2\setminus\{x\}= F_i$, and an element $y\in F_j\setminus (F_1\union F_2)$ such that $F_j=\{y\}\union (F_1\sect F_2)$. It follows that $F_i\union F_j=(F_1\union F_2)\union \{y\}$. This contradicts (ii).

(d)\implies (e): Assume that $|\Union_{i=1}^rF_i|=c+1$. After a suitable permutation of the elements of $[n]$  we may assume that $\Union_{i=1}^rF_i=\{1,\ldots,c+1\}$. Since $|F_i|=c$, there exists $j_i\in \{1,\ldots,c+1\}$ such that $F_i=\{1,\ldots,c+1\}\setminus\{j_i\}$. Since the sets $F_i$ are pairwise distinct it follows that $j_i\neq j_k$ for $i\neq k$. Thus after applying again suitable permutation we may assume that $j_i=i$ for $i=1,\ldots,r$.

The second statement follows similarly.

(e) \implies (f): In the first case, $v=2$ and $h=(c+1)$, while in the second case,
$v=r$ and $h=c-1+r$. Thus in both cases $\size I=n-c=\dim S/I$.

(f)\implies (g): By Lemma \ref{Lyubeznik} and the remark preceding the lemma, we have
\[
\depth S/L\geq \size L=\size I=\dim S/I=\dim S/L.
\]
Hence $S/L$ is Cohen-Macaulay.

Finally the implication (g) \implies (a) is trivial.
\end{proof}

\begin{Corollary}
\label{gorenstein}
With notation as above, the following conditions are equivalent:
\begin{enumerate}
\item[(a)] $J$ is a Gorenstein ideal for all choices of the  integers $a_{ij}$;
\item[(b)] $r=1$ or $c=1$.
\end{enumerate}
\end{Corollary}

\begin{proof} If $r=1$ or $c=1$, then $J$ is complete intersection for all choices of the integers $a_{ij}$. Thus (b) implies (a).

Conversely suppose condition (b) is not satisfied. We assume that $c>1$, and have to show that $r=1$.
By Theorem \ref{converse} we have $|\Sect_{i=1}^r F_i|=c-1$ or $|\Union_{i=1}^rF_i|=c+1$.

In the first case we may assume that $F_i=\{1,\ldots, c-1,i+c-1\}$ for $i=1,\ldots r$. Assume $r>1$, and let
$Q_{F_1}=(x_1^2,x_2,\ldots, x_c)$ and $Q_{F_i}=P_{F_i}$ for $i\geq 2$. Then $J=\Sect_{i=1}^rQ_{F_i}=(x_1^2,x_1x_2, \prod_{i=0}^{r-1}x_{c+i})$ is not Gorenstein, a contradiction.

In the second case suppose that $r\geq 3$. With the same argument as in the proof of Theorem \ref{converse} it follows that $I_A=\Sect_{i\in A}P_{F_i}$ is a Gorenstein ideal for all subsets $A\subset [r]$.
Therefore $P_{F_1}\sect P_{F_2}\sect P_{F_3}$ is Gorenstein. We may assume that
$F_1=\{1,2,\ldots, c\}$, $F_2=\{2,3,\ldots, c+1\}$ and $F_3=\{1,3,4,\ldots,c+1\}$.  Then $P_{F_1}\sect P_{F_2}\sect P_{F_3}=(x_1x_2, x_1x_{c+1}, x_2x_{c+1},x_3,\ldots,x_{c})$ is not Gorenstein, a contradiction.

On the other hand, if $r=2$, then $|\Sect_{i=1}^r F_i|=c-1$, and we are again in the first case. Thus we must have that $r=1$.
\end{proof}

\begin{Remark}
{\em From a  view point of Stanley-Reisner rings, the ideal $I$ in
the first case of condition (e) in Theorem \ref{converse}
corresponds to an iterated cone of a 0-dimensional simplicial
complex. In this case it is known that $S/I$ itself is Gorenstein if
the corresponding 0-dimensional simplicial complex consists of at
most 2 points, see \cite[Theorem 5.1(e)]{St}. The corollary also
follows from this fact.}
\end{Remark}


\newpage

\begin{thebibliography}{99}

\bibitem{BH} W.\ Bruns and J.\ Herzog, ``Cohen-Macaulay rings" (Revised edition), Cambridge Studies in  Advanced Mathematics {\bf 39}, Cambridge University Press, 1998.

\bibitem{BH1} W.\ Bruns and J.\ Herzog, On multigraded resolutions, Math.\ Proc.\ Camb.\ Phil.\ Soc.\ {\bf 118} (1995), 245--257.


\bibitem{BW}  A.\ Bj\"orner and M.\L.\ Wachs,  Shellable non-pure complexes and posets II, Trans.\ AMS {\bf 349} (1997) 3945--3975.


\bibitem{D} A.\ Dress, A new algebraic criterion for shellability, Beitr\"age zur Algebra und Geometrie {\bf 34} (1993), 45--55.

\bibitem{H} R.\ Hartshorne, Complete intersections in characteristic $p>0$, Amer.\ J.\ Math.\ {\bf 101} (1979), 380--383.

\bibitem{HS} J.\ Herzog and E.\ Sbarra, Sequentially Cohen-Macaulay modules and local cohomology, ``Algebra, arithmetic and geometry, Part I, II" (Mumbai, 2000), 327--340, Tata Inst. Fund. Res. Stud. Math., {\bf 16}, Tata Inst. Fund. Res., Bombay, 2002.

\bibitem{L} G.\ Lyubeznik, On the arithmetic rank of monomial ideals, J.\ Alg.\ {\bf 112} (1988), 86--89.

\bibitem{Sch} P.\ Schenzel, On the number of faces of simplicial complexes and the purity of Frobenius, Math.\ Z.\ {\bf 178} (1981), 125--142.


\bibitem{St} R.P.\ Stanley, ``Combinatorics and commutative algebra", Birkh\"auser, second edition, 1996.


\bibitem{T} Y.\ Takayama, A generalized Hochster's formula for local
cohomologies of monomial ideals, preprint 2004.

\bibitem{Ta} A.\ Taylor,   The inverse Groebner basis problem in
codimension two, J.\   Symb.\ Comp.\ {\bf 33} (2002), 221--238.

\bibitem{Te} N.\ Terai,  Local cohomology modules with respect to monomial ideals, preprint 1998.



\end{thebibliography}
\end{document}